\documentclass[a4paper,10pt,reqno, english]{amsart}

\usepackage[utf8]{inputenc}
\usepackage[T1]{fontenc}
\usepackage{amsmath,amsthm}
\usepackage{amsfonts,amssymb,enumerate}
\usepackage{url,paralist}
\usepackage{listings}
\usepackage{mathtools,color}
\usepackage[colorlinks=true,urlcolor=blue,linkcolor=red,citecolor=magenta,hypertexnames=false]{hyperref}
\usepackage{tikz}
\tikzset{auto}
\usepackage{pgfplotstable}

\usepackage{enumerate}
\usepackage{anysize}

\theoremstyle{plain}
\newtheorem{theorem}{Theorem}[section]
\newtheorem{lemma}[theorem]{Lemma}

\newtheorem{proposition}[theorem]{Proposition}

\newtheorem*{theorem*}{Theorem}

\newtheorem*{claim*}{Claim}
\newtheorem*{lemma*}{Lemma}

\theoremstyle{definition}
\newtheorem{definition}[theorem]{Definition}

\newtheorem{problem}[theorem]{Problem}

\DeclareMathOperator{\sign}{sign}
\DeclareMathOperator{\aff}{aff}
\newcommand{\conv}{\operatorname{conv}}
\newcommand{\defn}[1]{\emph{#1}} %
\newcommand\restr[2]{{ 
    \left.\kern-\nulldelimiterspace 
    #1 
    \vphantom{\big|} 
    \right|_{#2} 
}}

\definecolor{mygreen}{rgb}{0,0.6,0}
\definecolor{mygray}{rgb}{0.5,0.5,0.5}
\definecolor{mymauve}{rgb}{0.58,0,0.82}
\newcolumntype{R}[2]{%
    >{\adjustbox{angle=#1,lap=\width-(#2)}\bgroup}%
    l%
    <{\egroup}%
}
\begin{document}
\title[Ten colored points]{
    A new $k$-partite graph $k$-clique iterator and the optimal colored Tverberg problem for ten colored points}

\author[J.~Kliem]{Jonathan Kliem}
\address[J.~Kliem]{Institut f\" ur Mathematik, Freie Universität Berlin, Germany}
\email{jonathan.kliem@fu-berlin.de}

\thanks{
  J.K.~receives funding by the Deutsche Forschungsgemeinschaft DFG under Germany's Excellence Strategy – The Berlin Mathematics Research Center MATH+ (EXC-2046/1, project ID: 390685689).
}
\date{\today}

\begin{abstract}
    We provide an algorithm that verifies the optimal colored Tverberg problem for $10$ points in the plane:
    Every $10$ points in the plane in color classes of size at most $3$ can be partitioned in $4$ rainbow pieces such that their convex hulls intersect in a common point.

    This is achieved by translating the problem to $k$-partite graphs and using a new algorithm to verify that those graphs do not have a $k$-clique.
\end{abstract}

\maketitle

\tableofcontents


\section{Introduction}

Tverberg's Theorem has many variations.
This is a compact presentation of some of them:

\begin{problem}[Topological/Affine Tverberg]
    Let $d \geq 1$ and $r \geq 2$ be integers and let $N = (d+1)(r-1)$ and let $f \colon \Delta_N \to \mathbb{R}^d$ be a continous/affine map from the $N$-simplex to the $d$-dimensional Euclidean space.
    Are there pairwise disjoing faces $\sigma_1,\dots,\sigma_r$ of $\Delta_N$ such that their images with respect to $f$ intersect?

    Optimal Colored: Let further $C_1,\dots,C_m$ be color classes of the vertex set of size at most $r-1$. Can $\sigma_1,\dots,\sigma_r$ be chosen rainbow (for each $i=1,\dots,r$ vertices of $\sigma_i$ have pairwise distinct colors)?
\end{problem}

The affine Tverberg Problem without colors is solved affirmatively by Tverberg's theorem~\cite{Tverberg1966}.
For the topological Tverberg Problem the answer is more complicated and depends on $r$:
Bárany, Shlosman and Szücs provided in 1981 that the answer is yes when $r$ is a prime~\cite{Shlosman}.
In 1987, Özaydin extended this for $r$ a prime power~\cite{PrimePowers} in a never published preprint.
This first published proof was later provided by Volovikov~\cite{Volovikov}.
Recently, Blagojević, Frick and Ziegler discoverd that the answer is no for all $r$ that are not prime powers~\cite{BFZ2019}.

One generalization of the Tverberg Problem is the optimal colored version, which implies other previous generalizations.
Blagojević, Matschke and Ziegler~\cite{BMZ2015} showed that the topological optimal colored Tverberg problem holds for primes $r$.
It is unkown whether the optimal colored version holds for prime powers even in the affine case.

We refer the reader to Bárány, Blagojević and Ziegler~\cite{History} for the history of the problem.
In particular, \cite[Figure~10]{History} motivates this paper.
We show:

\begin{theorem}\label{Thm:TenPoints}
    Let $X = \{x_1,\dots,x_{10}\}$ be points in $\mathbb{R}^2$.
    Let $C_1,\dots,C_m$ be a partition of $X$, such that $|C_j| \leq 3$ for any $j = 1,\dots,m$.
    There is a partition $X_1,\dots,X_4$ of $X$ such that
    $|X_i \cap C_j| = 1$ for any $i = 1,\dots,4$, $j = 1,\dots,m$ and such that
    \[
        \conv(X_1) \cap \conv(X_2) \cap \conv(X_3) \cap \conv(X_4) \neq \emptyset.
    \]
\end{theorem}

We obtain this on the level of oriented matroids using a $k$-clique iterator on $k$-partite graphs.
As other algorithms fail to analyse those graphs, we present a new $k$-clique iterator for $k$-partite graphs.

As we will see in Proposition~\ref{Prop:StrongGeneralPosition}, we may assume the points $x_1,\dots,x_{10}$ to be in strong general position,
as moving them slightly the new Tverberg partitions will be a subset of the old ones.
We may also assume that $|C_1| \geq |C_2| \geq \dots \geq |C_m|$ and that $|C_m \cup C_{m-1}| > 3$.
So $(|C_1|,\dots, |C_m|)$ is one of $(3,3,3,1)$, $(3,3,2,2)$, $(2,2,2,2,2)$.
The cases
\[
    d = 2,\, r=4,\, \left(|C_1|,\dots,|C_m|\right) = (3,3,2,2) \quad \text{and} \quad
    d = 2,\, r=4,\, \left(|C_1|,\dots,|C_m|\right) = (2,2,2,2,2)
\]
can be reduced to
\[
    d = 3,\, r=4,\, \left(|C_1|,\dots,|C_m|\right) = (3,3,3,3,1) \quad \text{and} \quad
    d = 4,\, r=4,\, \left(|C_1|,\dots,|C_m|\right) = (3,3,3,3,3,1)
\]
by~\cite[Reduction of Thm.~2.2 to Thm~2.1]{BMZ2015}.
However, the higher dimensional cases with more points seems much harder to solve.

There are $2800 + 6300 + 945 = 10045$ partitions of $x_1,\dots,x_{10}$ into color classes of cardinalities $(3,3,3,1)$, $(3,3,2,2)$ and $(2,2,2,2,2)$, respectively.
So for each collection of points $x_1,\dots,x_{10} \in \mathbb{R}^2$ we must check for all those $10045$ color partitions, whether there is a rainbow Tverberg partition.

Tverberg partitions of $10$ points in general position in $\mathbb{R}^d$ are either of cardinalities $(3,3,3,1)$ or $(3,3,2,2)$.
We will call those types $3,3,3,1$ and $3,3,2,2$.
For type $3,3,3,1$ we can simply check the oriented matroid on $x_1,\dots,x_{10}$ induced by the point configuration.
For type $3,3,2,2$, this does not suffice: All $10$-gons have the same oriented matroid, but their Tverberg partitions are different.

This paper is organized as follows:
In Section~\ref{Sec:StrongGeneralPosition} we reduce the problem to points in strong general position.
In Section~\ref{Sec:AbstractGraph} we will explain how Theorem~\ref{Thm:TenPoints} reduces to verifying for each acyclic chirotope on $x_1,\dots,x_{10}$ of rank $3$ that some $k$-partite graph does not have a $k$-clique.
In Section~\ref{Sec:ConcreteGraph} we will show how a slightly modified graph can be efficiently constructed.
Finally, in Section~\ref{Sec:kpkc} we will introduce a new $k$-partite $k$-clique iterator that is able to verify this.
We compare its performance with other algorithms on random graphs in Sectoin~\ref{Sec:Benchmarks}.

\section{Reducing the problem to strong general position}\label{Sec:StrongGeneralPosition}

To simplify the combinatorics, it would be nice if the points $x_1,\dots,x_n$ in $\mathbb{R}^d$ are in general position:
\begin{definition}[{\cite[Def.~2.1]{Perles}}]
    A finite set $X = \{x_1,\dots,x_n\}$ in $\mathbb{R}^d$ is said to be in \defn{strong general position}
    if every subset of $X$ of size $\leq d+1$ is affinely independent and
    for any collection $\{X_1,\dots,X_r\}$ of $r$ pairwise disjoint subsets of $X$ we have
    \[
        d - \dim \bigcap_{i=1}^r \aff(F_i) = \min\left(d+1, \sum_{i=1}^r (d - \dim \aff(F_i)\right).
    \]
\end{definition}
In the case of points in $\mathbb{R}^2$, this means that no three lines (defined by pairwise distinct points of $X$) contain a common point.

\smallskip

We will show that we can slightly modify our points without creating new Tverberg partitions
to achieve this property:

\begin{lemma}\label{lem:dense}
    The subset of points in strong general position of $(\mathbb{R}^d)^n$ is dense.
\end{lemma}
\begin{proof}
    Perles and Sigron state this implicitly in~\cite[Sec.~3]{Perles} by showing that the points in strong general position can be seen as the subset on which a
    non-trivial polynomial in all points evaluates non-zero.
\end{proof}

\begin{lemma}\label{lem:open}
    Let $I_1,\dots,I_r$ be disjoint subsets of $\{1,\dots,n\}$.
    The set $U_{I_1,\dots,I_r}$ of points $(x_1,\dots,x_n)$ in $(\mathbb{R}^d)^n$ for which
    \[
        \emptyset = \bigcap_{i=1}^r \conv \left\{ x_j \colon j \in I_i\right\}
    \]
    is open.
\end{lemma}
\begin{proof}
    Let $D$ be the minimal distance that can be achieved from the convex sets induced by $I_1,\dots, I_r$:
    \[
        D = \min_{z \in \mathbb{R}^d} \max_{i=1,\dots,r} \Big( d\left(z, \conv \left\{ x_j \colon j \in I_i\right\} \right) \Big).
    \]
    This distance $D$ being non-zero is equivalent to
    \[
        \emptyset = \bigcap_{i=1}^r \conv \left\{ x_j \colon j \in I_i\right\}.
    \]
    Suppose that $(x_1,\dots,x_n) \in U_{I_1,\dots,I_r}$, which implies $D > 0$.
    Let $(y_1,\dots,y_r) \in (\mathbb{R}^d)^n$ with $d(x_i,y_i) < D$ for all $i = 1,\dots,n$.
    Then
    \begin{align*}
        & \min_{z \in \mathbb{R}^d} \max_{i=1,\dots,r} \Big( d\left(z, \conv \left\{ y_j \colon j \in I_i\right\} \right) \Big) \\
        > & \min_{z \in \mathbb{R}^d} \max_{i=1,\dots,r} \Big( d\left(z, \conv \left\{ x_j \colon j \in I_i\right\} \right) \Big) - D\\
        \geq & 0.
    \end{align*}
    This implies that
    \[
        \emptyset = \bigcap_{i=1}^r \conv \left\{ y_j \colon j \in I_i\right\}.
    \]
    and therefore $(y_1,\dots,y_r) \in U_{I_1,\dots,U_r}$.
\end{proof}

\begin{proposition}\label{Prop:StrongGeneralPosition}
    Let $X = \{x_1,\dots,x_n\}$ be a set of points in $\mathbb{R}^d$.
    There exists a set $Y = \{y_1,\dots,y_n\}$ in $\mathbb{R}^d$ in strong general position such that
    for any disjoint subsets $I_1,\dots,I_r$ of $\{1,\dots,10\}$
    with
    \[
        \emptyset = \bigcap_{i=1}^r \conv \left\{ x_j \colon j \in I_i\right\}.
    \]
    we have
    \[
        \emptyset = \bigcap_{i=1}^r \conv \left\{ y_j \colon j \in I_i\right\}.
    \]
\end{proposition}
\begin{proof}
    There exists only finitely many disjoint subsets $I_1,\dots,I_r$ with
    \[
        \emptyset = \bigcap_{i=1}^r \conv \left( x_j \colon j \in I_i\right).
    \]
    For those the sets $U_{I_1,\dots,I_r}$ are open by Lemma~\ref{lem:open} and their intersection contains $(x_1,\dots,x_n)$.
    This open non-empty intersection must contain some $(y_1,\dots,y_n)$ in strong general position by Lemma~\ref{lem:dense}.
\end{proof}

\section{Reducing the problem to $k$-partite graphs}\label{Sec:AbstractGraph}

To utilize a computer, we use the chirotope axioms of an oriented matroid.
For an introduction on oriented matroids we refer the reader to Björner, Las Vergnas, Sturmfels, White and Ziegler~\cite{OM}.

\begin{definition}[{\cite[Def. 3.5.3]{OM}}]
    A \defn{chirotope} of rank $r$ on a set $E$ is a mapping $\chi \colon E^r \to \{-1,0,1\}$, which satisfies the following three properties:
    \begin{enumerate}[({B}1)]
            \setcounter{enumi}{-1}
        \item $\chi$ is not identically zero,
        \item $\chi$ is alternating, that is
            \[
                \chi(x_{\sigma_1},\dots,x_{\sigma_r}) = \sign(\sigma) \chi(x_1,\dots,x_r)
            \]
            for all $x_1,\dots,x_r \in E$ and every permutation $\sigma$,
        \item for all $x_1,\dots,x_r, y_1,\dots,y_r \in E$ such that
            \[
                \chi(y_i,x_2,x_3,\dots,x_r) \cdot \chi(y_1,y_2,\dots,y_{i-1},x_1, y_{i+1},y_{i+2},\dots,y_r) \geq 0
            \]
            for $i=1,\dots,r$, we have
            \[
                \chi(x_1,\dots,x_r) \cdot \chi(y_1,\dots,y_r) \geq 0.
            \]
    \end{enumerate}
\end{definition}

Instead of Axiom (B2) will will use an equivalent formulation:

\begin{lemma}[{\cite[Lem~3.5.4]{OM}}]\label{Lem:B2}
    Let $\chi \colon E^r \to \{-1,0,1\}$ be a map satisfying (B0) and (B1).
    Then (B2) is equivalent to the following:
    For any $\chi(x_1,\dots,x_r)\chi(y_1,\dots,y_r) \neq 0$ there exists $i \in \{1,\dots,r\}$ such that
    \[
        \chi(x_1,\dots,x_r)\chi(y_1,\dots,y_r) = \chi(y_i,x_2,x_3,\dots,x_r)\chi(y_1,\dots,y_{i-1},x_1,y_{i+1},\dots,y_r).
    \]
\end{lemma}

As explained on \cite[Page 5]{OM}, every point configuration in $\mathbb{R}^d$ corresponds to an acyclic chirotope of rank $d+1$.
We will provide a definition of an acyclic chirotope of rank $3$:

\begin{definition}
    Let $\chi$ be a chirotope on $E$ of rank $3$.
    The chirotope $\chi$ is \defn{acyclic} if for every $x_1,\dots,x_4 \in E$ with $\chi(x_1,x_2,x_3) \neq 0$
    one of
    \[
        \chi(x_1,x_2,x_4), \,\chi(x_2,x_3,x_4), \,\chi(x_3,x_1,x_4)
    \]
    is equal to $\chi(x_1,x_2,x_3)$.
\end{definition}
\begin{lemma}
    The definition of an acyclic chirotope of rank $3$ agrees with the definition~\cite[Def.~3.4.7]{OM}:
    An oriented matroid is acyclic if it does not contain a positive circuit.
\end{lemma}
\begin{proof}
    Let $C$ be a circuit with support contained in $x_1,\dots,x_4$.
    As $x_1,x_2,x_3$ is a basis, we may assume that $x_4 \in C^{+}$.

    $x_1 \in \underline{C}$ is equivalent to $x_2, x_3, x_4$ being a basis.

    So either $\chi(x_2, x_3, x_4) = 0$ or
    by \cite[Prop.~3.5.2]{OM} we conclude that
    \begin{align*}
        \chi(x_2, x_3, x_4) = \chi(x_4, x_2, x_3) = -C(x_1)C(x_4) = -C(x_1)\chi(x_1,x_2,x_3).
    \end{align*}

    Our statements have been invariant with respect to cyclic permutation of $x_1, x_2, x_3$.
    If, $C$ is a positive circuit, then all of
    \[
        \chi(x_1,x_2,x_4), \,\chi(x_2,x_3,x_4), \,\chi(x_3,x_1,x_4)
    \]
    are zero or equal to $-\chi(x_1, x_2, x_3)$.
    On the other hand, if $C$ is not a positive circuit, then one of them is equal to $\chi(x_1, x_2, x_3)$.
\end{proof}

Let $x_1,x_2,\dots,x_{10} \in \mathbb{R}^2$ in strong general position.
For any distinct $a,b,c,d \in \{x_1,\dots,x_{10}\}$ we denote the intersection of the lines $ab$ and $cd$ by $y_{a,b,c,d}$.
Note that $y_{a,b,c,d}$ has $8$ different notations:
\[
    y_{a,b,c,d},\, y_{a,b,d,c},\, y_{b,a,c,d},\, y_{b,a,d,c},\, y_{c,d,a,b},\, y_{c,d,b,a},\, y_{d,c,a,b},\, y_{d,c,b,a}.
\]
We define the sets
\[
    X := \{x_1,\dots,x_{10}\}, \quad Y := \left \{y_{a,b,c,d} \colon |\{a,b,c,d\}| = 4, \{a,b,c,d\} \subset\{1,\dots,10\}\right\}.
\]

This way the points $x_1,\dots,x_{10} \in \mathbb{R}^2$ induce an acyclic chirotope of rank $3$ on $X \sqcup Y$ with
\[
    \chi(x, x', x'') \neq 0
\]
for all $x, x', x'' \in X$ and
\[
    \chi(a, b, y_{a,b,c,d}) = 0 = \chi(c,d, y_{a,b,c,d})
\]
for all $a,b,c,d \in \{x_1,\dots,x_{10}\}$ pairwise distinct.

\begin{definition}
    The point $y_{a,b,c,d}$ is an \defn{intersection point}, if $y_{a,b,c,d} = \conv(a,b) \cap \conv(c,d)$.
\end{definition}
\begin{lemma}\label{Lem:IntersectionPoint}
    The point $y_{a,b,c,d}$ is an intersection point, if and only if
    \[
        \chi(a, b, c)\chi(a,b, d) = -1 = \chi(c, d, a)\chi(c, d, b).
    \]
\end{lemma}
\begin{proof}
    $y_{a,b,c,d}$ is defined to be the intersection of the lines $ab$ and $cd$.
    Thus, it is an intersection point, if and only if $\conv(a,b) \cap \conv(c,d) \neq \emptyset$.
    As $a,b,c,d$ are in general position there is a cycle $C$ on the ground set $(a, b, c, d)$.
    By \cite[Prop.~3.5.2]{OM} we conclude that
    \begin{align*}
        \chi(a, b, c) &= -C(c)C(d)\chi(a, b, d)\\
        \chi(c, d, a) &= -C(a)C(b)\chi(c, d, b).
    \end{align*}
    Hence $C = (a, b, -c, -d)$ or $C = (-a, -b, c, d)$ is equivalent to
    \[
        \chi(a, b, c)\chi(a,b, d) = -1 = \chi(c, d, a)\chi(c, d, b).
    \]
\end{proof}

As $x_1,\dots,x_{10}$ are assumed to be in strong general position,
we have $\chi(y_{a,b,c,d}, e, f) \neq 0$ unless $\{e, f\} = \{a, b\}$ or $\{e, f\} = \{c, d\}$.

A chirotope on $X \sqcup Y$ determines all Tverberg partitions:
\begin{lemma}\label{Lem:Restriction}
    Let $a,b,c,d,e,f,g \in \{x_1,\dots,x_{10}\}$ be pairwise distinct and let $y_{a,b,c,d}$ be an intersection point.
    \begin{enumerate}
        \item\label{Restriction1}
            $a \in \conv(e, f, g)$ is equivialent to
            \[
                \chi(e, f, a) = \chi(f, g, a) = \chi(g, e, a).
            \]
        \item\label{Restriction2}
            \[
                y_{a,b,c,d} \in \conv(e,f, g)
            \]
            is equivalent to
            \[
                \chi(e, f, y_{a,b,c,d}) = \chi(f, g, y_{a,b,c,d}) = \chi(g, e, y_{a,b,c,d}).
            \]
    \end{enumerate}
\end{lemma}
\begin{proof}
    $a,e,f,g$ and $y_{a,b,c,d}, e, f, g$ are in general position by assumptions.
    Hence,~\eqref{Restriction2} is implied by~\eqref{Restriction1}.
    As $a, e, f, g$ are in general position, there exists a cycle $C$ with ground set $a, e, f, g$.
    W.l.o.g. we have $C(g) = 1$.
    By \cite[Prop.~3.5.2]{OM} we conclude that
    \begin{align*}
        \chi(e, f, a) = C(f)C(g)\chi(g, e, a) = C(f)\chi(g, e, a)\\
        \chi(e, f, a) = C(e)C(g)\chi(f, g, a) = C(e)\chi(f, g, a)
    \end{align*}

    $a \in \conv(e, f, g)$ is equivalent to $C = (-a, e, f, g)$, which implies
    \[
        \chi(e, f, a) = \chi(f, g, a) = \chi(g, e, a).
    \]
    For the other direction, this equality implies $C(e) = 1 = C(f)$.
    As the oriented matroid is acyclic we have $C(a) = -1$.
\end{proof}

This means that restriction of the chirotope to $X$ determines, whether there is a Tverberg partition of type $3,3,3,1$ and the restriction
to each $X \cup \{y_{a,b,c,d}\}$ determines whether there is a Tverberg partition of type $3,3,2,2$.

The restrictions of those chirotopes to $X$ are exactly acyclic chirotopes on $x_1,\dots,x_{10}$ of rank $3$.
They were previously classified by Aichholzer, Aurenhammer and Krasser~\cite{Aichholzer2001}.
In total there are 14,320,182 such chirotopes and 14,309,547 of them are realizable.

By Lemma~\ref{Lem:Restriction} the chirotope $\chi$ determines all Tverberg partitions.
The restriction of $\chi$ to $X$, determines all partitions of type $3,3,3,1$, but partitions of type $3,3,2,2$ are not determined other than in special cases.
We give a definition to allow us to talk about this:

\begin{definition}
    Let $X \subseteq Z \subseteq X \sqcup Y$.
    The Tverberg partitions on $\restr{\chi}{Z}$ are those Tverberg partitions that exist by Lemma~\ref{Lem:Restriction}:
    \begin{itemize}
        \item Tverberg partitions on $\restr{\chi}{Z}$ contain all partitions of type $3,3,3,1$ regardless of the choice of $Z$.
        \item Tverberg partitions on $\restr{\chi}{Z}$ contain exactly those partitions $(X_1,X_2, \{a, b\}, \{c,d\})$ of type $3,3,2,2$,
            for which $y_{a,b,c,d} \in Z$.
    \end{itemize}
\end{definition}

For the remainder $y := y_{a,b,c,d}$.
Depending on $\restr{\chi}{X}$, some of the intersection points can be ignored:
\begin{lemma}\label{Lem:ValidIntersectionPoint}
    Let $a,b,c,d,e,f,g \in \{x_1,\dots,x_{10}\}$ be pairwise distinct such that $y := y_{a,b,c,d}$ is an intersection point.
    \begin{enumerate}
        \item
            If $y \in \conv(e,f,g)$ then
            \[
                \chi(a, b, h) \geq 0
            \]
            for some $h \in \{e,f,g\}$.
        \item
            If there is a Tverberg partition of type $3,3,2,2$ with $y = y_{a,b,c,d}$,
            then
            \[
                \left | \left \{h \in \{x_1,\dots,x_{10}\} \setminus \{a,b,c,d\} \colon \, \chi(a,b,h) = 1 \right \} \right | \geq 2.
            \]
    \end{enumerate}
\end{lemma}
\begin{proof}
    \begin{enumerate}
        \item
            Suppose that
            \[
                \chi(a,b,f) = \chi(a,b,g) = -1.
            \]
            Then, $y \in \conv(e,f,g)$ implies by Lemma~\ref{Lem:Restriction} that
            \[
                \chi(e, f, y) = \chi(f, g, y) = \chi(g, e, y).
            \]
            We conclude
            \begin{align*}
                \chi(y, a, b) \chi(e, f, g) &= 0,\\
                \chi(f, a, b) \chi(y, e, g) & \\
                = \chi(g, a, b) \chi(y, f, e) &\neq 0.
            \end{align*}
            With Lemma~\ref{Lem:B2} and $\chi(e,a,b)\chi(y,f,g) \neq 0$ we conclude that
            \[
                \chi(e, a, b) \chi(y, f, g) = \chi(f, a, b) \chi(y, e, g).
            \]
            or
            \[
                \chi(e, a, b) = -\chi(f, a, b).
            \]
        \item Follows.
    \end{enumerate}
\end{proof}
This means that Tverberg partitions of type $3,3,2,2$ with $y = y_{a,b,c,d}$ can only exist,
if at least $2$ of the remaining points lie above the line $ab$.
By symmetry also $2$ points must lie below and the the same applies for the line $cd$.
We summarize this in the following definition:
\begin{definition}
    An intersection point $y = y_{a,b,c,d}$ is \defn{valid} if
    \[
        2 \leq \left | \left \{h \in \{x_1,\dots,x_{10}\} \setminus \{a,b,c,d\} \colon \, \chi(a,b,h) = 1 \right \} \right | \leq 4
    \]
    and
    \[
        2 \leq \left | \left \{h \in \{x_1,\dots,x_{10}\} \setminus \{a,b,c,d\} \colon \, \chi(c,d,h) = 1 \right \} \right | \leq 4.
    \]
\end{definition}

According to Lemma~\ref{Lem:ValidIntersectionPoint} we only need to consider valid intersection points.
Depending on $\restr{\chi}{X}$ there are up to $70$ valid intersection points.
This maximum is attained e.g.~for $\restr{\chi}{X}$ being the chirotope of a $10$-gon.
\begin{definition}
    Let $\restr{\chi}{X}$ be a chirotope on $X$.
    We construct a graph $G'(\restr{\chi}{X})$ as follows.
    For each valid intersection point $y_{a,b,c,d}$ we have vertices for each extension of $\restr{\chi}{X}$ to $\restr{\chi}{X \cup \{y_{a,b,c,d}\}}$.
    There is an edge between two such extensions, $\restr{\chi}{X \cup \{y_{a,b,c,d}\}}$ and $\restr{\chi}{X \cup \{y_{a',b',c',d'}\}}$
    if they are restrictions of a chirotope on $X \cup \{y_{a,b,c,d}, y_{a',b',c',d'}\}$.
\end{definition}
This graph is a $(k-1)$-partite graph, where $k-1$ is the number of valid intersection points.

Any chirotope $\chi$ of rank $3$ on $X \sqcup Y$ corresponds to a $(k-1)$-clique in $G'(\restr{\chi}{X})$.
According to Lemma~\ref{Lem:Restriction} all the Tverberg partitions are determined by restrictions of the chirotope corresponding to the $(k-1)$-clique.

\begin{definition}
    Let $\mathcal{C}$ be the set of all color partitions of cardinalities $(3,3,3,1)$, $(3,3,2,2)$, $(2,2,2,2)$.
    We construct a graph $G(\restr{\chi}{X})$ from $G'(\restr{\chi}{X})$ as follows:
    We add $\mathcal{C}$ to the vertices.
    We add the following edges:
    For each $(C_1,\dots,C_m) \in \mathcal{C}$ and each $\restr{\chi}{X \cup \{y_{a,b,c,d}\}}$
    we add an edge if {\bf none} of the Tverberg partitions on $\restr{\chi}{X \cup \{y_{a,b,c,d}\}}$ is rainbow with respect to $(C_1,\dots,C_m)$.
\end{definition}

This graph is $k$-partite, where $k-1$ is the number of valid intersection points.
In section~\ref{Sec:kpkc} we will develop an algorithm, which shows that none of the graphs has a $k$-clique:

\begin{proposition}\label{Prop:KClique}
    For each $\restr{\chi}{X}$ the graph $G(\restr{\chi}{X})$ does not have a $k$-clique, where $k-1$ is the number of valid intersection points.
\end{proposition}
This proposition shows our main Theorem:
\begin{proof}[{Proof of Theorem~\ref{Thm:TenPoints}}]
    Let $x_1,\dots,x_{10}$ in $\mathbb{R}^2$ with induced chirotope $\chi$ on $X \sqcup Y$
    and let $(C_1,\dots,C_m) \in \mathcal{C}$ be a color partition.
    Let $y_1,\dots,y_{k-1}$ be the valid intersection points of $\restr{\chi}{X}$.
    The restrictions to the valid intersection points $\restr{\chi}{X \cup \{y_1\}}, \dots, \restr{\chi}{X \cup \{y_{k-1}\}}$ form a $(k-1)$-clique of $G'(\restr{\chi}{X})$.

    However, by Proposition~\ref{Prop:KClique} this $(k-1)$-clique does not extend to a $k$-clique of $G(\restr{\chi}{X})$.
    In particular $(C_1,\dots,C_m)$ is not connected to some $\restr{\chi}{X \cup \{y_{a,b,c,d}\}}$.
    This means that a Tverberg partition on $\restr{\chi}{X \cup \{y_{a,b,c,d}\}}$ is rainbow with respect to $(C_1,\dots,C_m)$.
\end{proof}

\section{Obtaining the graph}\label{Sec:ConcreteGraph}

Instead of constructing $G(\restr{\chi}{X})$ we will construct a slightly larger graph $H(\restr{\chi}{X})$, which is simpler to compute:

\begin{lemma}\label{Lem:EnlargedGraph}
    Let $G$ be a graph. If $G$ has a $k$-clique, then so has any graph that is obtained by
    \begin{itemize}
        \item adding vertices,
        \item adding edges,
        \item vertex identification of pairwise non-adjacent vertices.
    \end{itemize}
\end{lemma}
\begin{proof}
    Trivial.
\end{proof}

We do not check all axioms when determining all $\restr{\chi}{X \cup \{y_{a,b,c,d}\}}$.
Thus we might add some vertices.
When checking whether $\restr{\chi}{X \cup \{y_{a,b,c,d}\}}$ and $\restr{\chi}{X \cup \{y_{a',b',c',d'}\}}$ have a common extension, we do not verify all axioms.
Thus, we might add some edges.
Instead of having a vertex for each $\restr{\chi}{X \cup \{y_{a,b,c,d}\}}$, we will only determine $\chi$ for some values that suffice to determine all Tverberg partitions.
This operation might identify vertices.
All identified vertices are in the same part and are therefore pairwise non-adjacent.

\subsection{The vertices of the graph}

For the remainder, $a,b,c,d,e,f,g \in \{x_1,\dots,x_{10}\}$ will be pairwise distinct such that $y := y_{a,b,c,d}$ is a valid intersection point.
In addition, $\{i,j,k\}$ is some permutation of $\{e, f, g\}$.

\begin{figure}
    \begin{center}
        \begin{tikzpicture}[scale=.7]
            \node (a) at (-2, 0) [circle, fill=black, scale=0.2, label=left:$a$] {};
            \node (b) at ( 2, 0) [circle, fill=black, scale=0.2, label=right:$b$] {};
            \node (c) at (0, -1) [circle, fill=black, scale=0.2, label=below:$c$] {};
            \node (d) at (0,  2) [circle, fill=black, scale=0.2, label=above:$d$] {};
            \node (e) at (-1, 0.9) [circle, fill=black, scale=0.2, label=above:$e$] {};
            \node (f) at ( 1, 1.2) [circle, fill=black, scale=0.2, label=above:$f$] {};
            \draw (a) -- (b);
            \draw (c) -- (d);
            \draw (e) -- (f);
        \end{tikzpicture}
    \end{center}
    \caption{Sketch of Lemma~\ref{Lem:SameSide}.}
    \label{Fig:SameSide}
\end{figure}
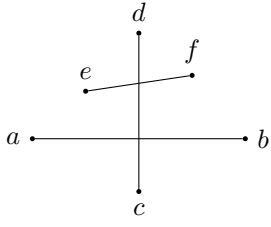

\begin{lemma}\label{Lem:PointInSegment}
    For any $h \in X \sqcup Y$ it holds that $\chi(a,b,h) = \chi(a,y,h) = \chi(y,b,h)$.
\end{lemma}
\begin{proof}
    As $y$ is a valid intersection point, we have that $\chi(a,c,d) \neq \chi(c,d,b)$.
    We inspect
    \begin{align*}
        \chi(y,c,d)\chi(a,b,g) &= 0,\\
       \chi(b,c,d)\chi(y,a,g) \\
        = \chi(g,c,d)\chi(y,b,a) &\neq 0.
    \end{align*}
    Lemma~\ref{Lem:B2} implies that
    \[
        \chi(a,c,d)\chi(y,b,g) = \chi(b,c,d)\chi(y,a,g).
    \]
    We conclude that $\chi(y,b,g) = \chi(a,y,g)$.
    If $\chi(a,b,g) \neq 0$, then one of
    \[
        0 = \chi(a,b,y), \quad \chi(b,g,y) = \chi(g,a,y)
    \]
    must be equal to $\chi(a,b,g)$ as the chirotope is acyclic.
    If $\chi(a,y,g) \neq 0$, then one of
    \[
        \chi(a,y,b) = 0, \quad \chi(y,g,b) = -\chi(a,y,g), \quad \chi(g,a,b)
    \]
    must be equal to $\chi(a,y,g)$.
\end{proof}
\begin{lemma}[{See Figure~\ref{Fig:SameSide}}]\label{Lem:SameSide}
    Suppose that $\chi(a,b,e) = \chi(a,b,f)$ and $\chi(c,d,e) \neq \chi(c,d,f)$ it follows that
    \[
        \chi(e, f, y) = -\chi(a,b,c)\chi(c,d,f)\chi(a,b,e) = \chi(c,d,a)\chi(c,d,f)\chi(a,b,e).
    \]
\end{lemma}
\begin{proof}
    By applying Lemma~\ref{Lem:PointInSegment} we have that
    \begin{align*}
        \chi(y, y, c)\chi(b, e, f) &= 0 \\
        \chi(e, y, c) \chi(y, b, f) &= \chi(e, d, c)\chi(a, b, f)\\
        = \chi(f, y, c) \chi(y, e, b) &= \chi(f, d, c)\chi(a, e, b) \neq 0.
    \end{align*}
    We conclude with Lemma~\ref{Lem:B2} and $\chi(b,y,c)\chi(y,e,f) \neq 0$ that
    \[
        \chi(b,y, c)\chi(y,e,f) = \chi(f, d, c)\chi(a,e,b).
    \]
    The statement follows with Lemma~\ref{Lem:PointInSegment}:
    \[
        \chi(b,y,c) = \chi(b,a,c) = \chi(y,a,c) = \chi(d,a,c).
    \]
\end{proof}

Depending on $\chi(a,b,e), \chi(c,d,e)$ the point $e$ is in one of four regions.
\begin{definition}
    If $\chi(a,b,e) = \chi(a,b,f)$ and $\chi(c,d,e) = \chi(c,d,f)$, then $e$ and $f$ are in the \defn{same region} with respect to $a,b,c,d$.
    If both equalities are false, then they are in \defn{opposite regions}.
    If exactly one equality holds, they are in \defn{neighboring regions}.
\end{definition}

\begin{lemma}\label{Lem:OppositeRegions}
    Suppose $y \in \conv(e,f,g)$.
    It follows that there exist $i,j \in \{e,f,g\}$
    such that $i$ and $j$ are in opposite regions with respect to $a,b,c,d$.
\end{lemma}
\begin{proof}
    Follows from Lemma~\ref{Lem:ValidIntersectionPoint} by pigeonhole principle.
\end{proof}

Let $i$ and $j$ be in opposite regions and let $y \in \conv(i,j,k)$.
Then $k$ is in the same region as $i$ or $j$ or in a neighboring region to both of them.

By relabeling, we will assume that $j$ and $k$ are in the same region or that
$\chi(a,b,i) = \chi(a,b,k)$ and $\chi(c,d,i) \neq \chi(c,d,k)$.

\begin{proposition}[{See Figure~\ref{Fig:OppositeRegions}}]\label{Prop:OppositeRegions} \
    \begin{enumerate}
        \item
            Suppose that $j$ and $k$ are in the same region and $i$ is in the opposite region to both of them with respect to $a,b,c,d$.
            The point $y$ is contained in $\conv(i,j,k)$ if and only if
            \[
                \chi(i,j,y) \neq \chi(i,k,y).
            \]
        \item
            Suppose that
            \[
                \chi(a,b,i) = \chi(a,b,k) \neq \chi(a,b,j), \quad
                \chi(c,d,i) \neq \chi(c,d,j) = \chi(c,d,k).
            \]
            The point $y$ is contained in $\conv(i,j,k)$ if and only if
            \[
                \chi(i, j, y) = -\chi(a,b,c)\chi(c,d,i)\chi(a,b,i).
            \]
    \end{enumerate}
\end{proposition}
\begin{proof}
    \begin{enumerate}
        \item
            If $y \in \conv(i,j,k)$ then certainly $\chi(i,j,y) \neq \chi(k,i,y)$ by Lemma~\ref{Lem:IntersectionPoint}.
            On the other hand, assume that $\chi(i,j,y) \neq \chi(i,k,y)$.
            We now have
            \begin{align*}
                \chi(j,c,d)\chi(i,k,y) & \\
                =\chi(k,c,d)\chi(j,i,y) &\neq 0,\\
                \chi(y, c, d)\chi(j,k,i) &= 0,
            \end{align*}
            which implies by Lemma~\ref{Lem:B2} and $\chi(i,c,d)\chi(j,k,y) \neq 0$ that
            \[
                \chi(i,c,d)\chi(j,k,y) = \chi(j,c,d)\chi(i,k,y).
            \]
            This yields $\chi(j,k,y) = -\chi(i,k,y)$ and therefore by
            Lemma~\ref{Lem:IntersectionPoint} we have that $y \in \conv(i,j,k)$.
        \item According to Lemma~\ref{Lem:SameSide} we have that
            \[
                \chi(j, k, y_{c,d,a,b}) = \chi(a,b,c)\chi(a,b,k)\chi(c,d,j)
            \]
            and that
            \[
                \chi(k, i, y) = -\chi(a,b,c)\chi(c,d,i)\chi(a,b,k).
            \]
            As $\chi(c,d,j) = -\chi(c,d,i)$ we conclude
            \[
                \chi(j,k,y) = \chi(k,i,y).
            \]
            Now according to Lemma~\ref{Lem:Restriction}, $y \in \conv(i,j,k)$ is equivalent to
            \[
                \chi(i,j,y) = \chi(j,k,y),
            \]
            which is by $\chi(a,b,i) = \chi(a,b,k)$ equivalent to
            \[
                \chi(i, j, y) = -\chi(a,b,c)\chi(c,d,i)\chi(a,b,i).
            \]
    \end{enumerate}
\end{proof}
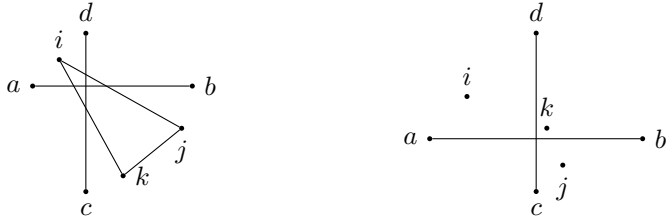
\begin{figure}
    \begin{center}
        \begin{tikzpicture}[scale=.7]
            \node (a) at (-1, 0) [circle, fill=black, scale=0.2, label=left:$a$] {};
            \node (b) at ( 2, 0) [circle, fill=black, scale=0.2, label=right:$b$] {};
            \node (c) at (0, -2) [circle, fill=black, scale=0.2, label=below:$c$] {};
            \node (d) at (0,  1) [circle, fill=black, scale=0.2, label=above:$d$] {};
            \node (i) at (-0.5, 0.5) [circle, fill=black, scale=0.2, label=above:$i$] {};
            \node (k) at (0.7, -1.7) [circle, fill=black, scale=0.2, label=right:$k$] {};
            \node (j) at (1.8, -0.8) [circle, fill=black, scale=0.2, label=below:$j$] {};
            \draw (a) -- (b);
            \draw (c) -- (d);
            \draw (i) -- (k) -- (j) -- (i);
        \end{tikzpicture}
        \hspace{2cm}
        \begin{tikzpicture}[scale=.7]
            \node (a) at (-2, 0) [circle, fill=black, scale=0.2, label=left:$a$] {};
            \node (b) at ( 2, 0) [circle, fill=black, scale=0.2, label=right:$b$] {};
            \node (c) at (0, -1) [circle, fill=black, scale=0.2, label=below:$c$] {};
            \node (d) at (0,  2) [circle, fill=black, scale=0.2, label=above:$d$] {};
            \node (i) at (-1.3, 0.8) [circle, fill=black, scale=0.2, label=above:$i$] {};
            \node (k) at ( 0.2, 0.2) [circle, fill=black, scale=0.2, label=above:$k$] {};
            \node (j) at (0.5, -0.5) [circle, fill=black, scale=0.2, label=below:$j$] {};
            \draw (a) -- (b);
            \draw (c) -- (d);
        \end{tikzpicture}
    \end{center}
    \caption{Sketch of Propositions~\ref{Prop:OppositeRegions} and~\ref{Prop:ExtensionObstruction}.}
    \label{Fig:OppositeRegions}
\end{figure}

We now have made precise what we mean by determining $\restr{\chi}{X \cup \{y\}}$ only for some values:
According to Proposition~\ref{Prop:OppositeRegions} it suffices to determine $\chi$ for all triples in $X$ and for all tripels $(i,j,y)$, where $i,j$ are in opposite regions with respect to $a,b,c,d$.

Now, we are not entirely free in choosing such an extension of $\chi$.
We will use some obstructions:

\begin{proposition}[{See Figure~\ref{Fig:OppositeRegions}}]\label{Prop:ExtensionObstruction} \
    \begin{enumerate}
        \item
            Suppose that $j$ and $k$ are in the same region and $i$ is in the opposite region to both of them with respect to $a,b,c,d$.
            \[
                \chi(i,j,y) \neq \chi(i,k,y)
            \]
            implies that
            \[
                \chi(i,j,y) = \chi(i,j,k).
            \]
        \item
            Suppose that
            \[
                \chi(a,b,i) = \chi(a,b,k) \neq \chi(a,b,j), \quad
                \chi(c,d,i) \neq \chi(c,d,j) = \chi(c,d,k).
            \]
            Now
            \[
                \chi(i, j, y) = -\chi(a,b,c)\chi(c,d,i)\chi(a,b,i).
            \]
            implies
            \[
                \chi(i, j, k) = -\chi(a,b,c)\chi(c,d,i)\chi(a,b,i)
            \]
    \end{enumerate}
\end{proposition}
\begin{proof}
    Both statements follow directly from Proposition~\ref{Prop:OppositeRegions} as the chirotope is acyclic.
\end{proof}

Also $\chi(i,j,y)$ is determined unless $\chi(i,j,a) \neq \chi(i,j,b)$ and $\chi(i,j,c) \neq \chi(i,j,d)$:

\begin{proposition}[{See Figure~\ref{Fig:NonCrossing}}]\label{Prop:NonCrossing}
    Let $i$ and $j$ be in opposite regions.
    The equality $\chi(i,j,a) = \chi(i,j,b)$ implies $\chi(i,j,a) = \chi(i, j, y)$.
\end{proposition}
\begin{proof}
    We have by assumptions and Lemma~\ref{Lem:PointInSegment}
    \begin{align*}
        \chi(a,i,j)\chi(y,b,i) &= \chi(a,i,j)\chi(a,b,i) \\
        = \chi(b,i,j)\chi(a,y,i) &= \chi(b,i,j)\chi(a,b,i) \neq 0,\\
        \chi(i,i,j)\chi(a,b,y) &= 0,
    \end{align*}
    which implies by Lemma~\ref{Lem:B2} and $\chi(y,i,j)\chi(a,b,i) \neq 0$ that
    \[
        \chi(y, i, j)\chi(a, b, i) = \chi(a, i, j)\chi(a, b, i).
    \]
\end{proof}
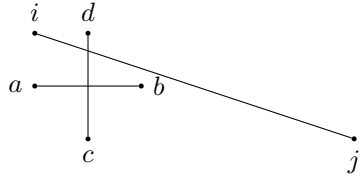
\begin{figure}
    \begin{center}
        \begin{tikzpicture}[scale=.7]
            \node (a) at (-1, 0) [circle, fill=black, scale=0.2, label=left:$a$] {};
            \node (b) at ( 1, 0) [circle, fill=black, scale=0.2, label=right:$b$] {};
            \node (c) at (0, -1) [circle, fill=black, scale=0.2, label=below:$c$] {};
            \node (d) at (0,  1) [circle, fill=black, scale=0.2, label=above:$d$] {};
            \node (i) at (-1, 1) [circle, fill=black, scale=0.2, label=above:$i$] {};
            \node (j) at (5, -1) [circle, fill=black, scale=0.2, label=below:$j$] {};
            \draw (a) -- (b);
            \draw (c) -- (d);
            \draw (i) -- (j);
        \end{tikzpicture}
    \end{center}
    \caption{Sketch of Proposition~\ref{Prop:NonCrossing}.}
    \label{Fig:NonCrossing}
\end{figure}

\begin{figure}
    \begin{center}
        \begin{tikzpicture}[scale=.7]
            \node (a) at (-3, 0) [circle, fill=black, scale=0.2, label=left:$a$] {};
            \node (b) at ( 1, 0) [circle, fill=black, scale=0.2, label=right:$b$] {};
            \node (c) at (0, -3) [circle, fill=black, scale=0.2, label=below:$c$] {};
            \node (d) at (0,  1) [circle, fill=black, scale=0.2, label=above:$d$] {};
            \node (i1) at (-1, 0.7) [circle, fill=black, scale=0.2, label=above:$i_1$] {};
            \node (j1) at (0.6, -1) [circle, fill=black, scale=0.2, label=right:$j_1$] {};
            \node (i2) at (-1.7, 0.4) [circle, fill=black, scale=0.2, label=above:$i_2$] {};
            \node (j2) at (0.3, -1.8) [circle, fill=black, scale=0.2, label=right:$j_2$] {};
            \draw (a) -- (b);
            \draw (c) -- (d);
            \draw (i1) -- (j1);
        \end{tikzpicture}
    \end{center}
    \caption{Sketch of Proposition~\ref{Prop:AdvancedExtensionObstruction}.}
    \label{Fig:AdvancedExtensionObstruction}
\end{figure}
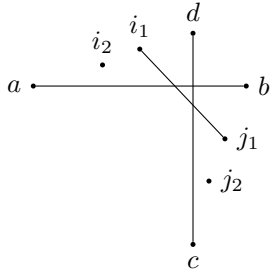

\begin{proposition}[{See Figure~\ref{Fig:AdvancedExtensionObstruction}}]\label{Prop:AdvancedExtensionObstruction}
    Suppose that $i_1,i_2$ are in the same region and $j_1,j_2$ are in the opposite region with respect to $a,b,c,d$.
    Also suppose that $\chi(i_1,j_1,a) \neq \chi(i_1,j_1,b)$ and $\chi(i_2, j_2, a) \neq \chi(i_2, j_2, b)$.
    Now
    \[
        \chi(i_1,j_1,i_2) = \chi(i_1,j_1,j_2) = -\chi(i_1,j_1,y)
    \]
    implies $\chi(i_2,j_2,y) = \chi(i_1,j_1,y)$.
\end{proposition}
\begin{proof}
    As the chirotope is acyclic one of $\chi(i_1,j_1,b), \chi(j_1,a,b), \chi(a,i_1,b)$ must be equal to $\chi(i_1,j_1,a)$.
    But $\chi(i_1,j_1,a) \neq \chi(i_1,j_1,b)$ and $\chi(j_1,a,b) = \chi(a,i_1,b)$.
    Hence it follows that $\chi(i_1,j_1,a) = \chi(a,b,j_1)$.
    Likewise $\chi(i_2,j_2,a) = \chi(a,b,j_2)$ and in particular $\chi(i_1,j_1, a) = \chi(i_2, j_2,a)$.

    We conclude that
    \begin{align*}
        \chi(y, i_1,j_1)\chi(a, i_2, j_2) &= \chi(i_1,j_1,y)\chi(a,b,j_2)\\
        = \chi(i_2, i_1,j_1)\chi(y, a, j_2) &= \chi(j_1,i_1,y)\chi(b, a, j_2)\\
        = \chi(j_2, i_1,j_1)\chi(y, i_2, a) &= \chi(j_1,i_1,y)\chi(b, i_2, a) \neq 0,\\
    \end{align*}
    which implies by Lemma~\ref{Lem:B2} and $\chi(a,i_1,j_1)\chi(y,i_2,j_2) \neq 0$ that
    \[
        \chi(a,i_1,j_1)\chi(y,i_2,j_2) = \chi(y,i_1,j_1)\chi(a, i_2,j_2).
    \]
\end{proof}

\subsection{The edges of the graph}

Given an acyclic chirotope $\restr{\chi}{X}$.
We are constructing a $k$-partite graph $H(\restr{\chi}{X})$ such that the existence of a $k$-clique in $G(\restr{\chi}{X})$ implies the existence of a $k$-clique in $H(\restr{\chi}{X})$.
This graph has the color partitions as vertices in the last part.
For each part corresponding to a valid intersection point $y$ it has a vertex $v$ for each collection of orientations of $(i,j,y)$ for $i$ and $j$ in opposite regions that agree with
Propositions~\ref{Prop:ExtensionObstruction},~\ref{Prop:NonCrossing}, and~\ref{Prop:AdvancedExtensionObstruction}.

Some vertices $v$ might not be extendable to a chirotope $\restr{\chi}{X \cup \{y\}}$.
Also it is possible that some vertices have multiple extensions.
However, $\restr{\chi}{X}$ determines the Tverberg partitions of types $(3,3,3,1)$ by Lemma~\ref{Lem:Restriction}.
Also the orientations fixed by $v$ determine all Tverberg partitions of type $(3,3,2,2)$ in $y$ by Proposition~\ref{Prop:OppositeRegions}.

Hence, $v$ determines the Tverberg partitions on $\restr{\chi}{X \cup \{y\}}$ and $v$ is connected to a color partition if none of those Tverberg partitions are rainbow with respect to this color partition.

Now let $w$ be another vertex corresponding to the valid intersection point $y' := y_{a',b',c',d'} \neq y$.
The vertices $v$ and $w$ must be connected by an edge, if the orientations are consistent.

We will always connect $v$ and $w$ unless $\{a,b\} = \{a', b'\}$ or $\{a, b\} = \{c', d'\}$ or $\{c,d\} = \{a', b'\}$ or $\{c, d\} = \{c',d'\}$.
By relabeling, those cases are reduced to $a = a'$, $b = b'$, $\chi(a,b,c) = \chi(a,b,c')$, which we will assume for the remainder of this section.
For each $v$ and $w$ we will check the following certificate:

\begin{proposition}[{See Figure~\ref{Fig:Edges}}]\label{Prop:Edges}
    Let $i, j$ be distinct from $a, b, c', d'$ such that
    $\chi(a,b,i) = \chi(a,b,c) = \chi(a,b,c')$ and $\chi(a,b,j) = \chi(a,b,d) = \chi(a,b,d')$.
    Suppose that $\chi(i,j,y) \neq \chi(c',d',y)$ and let $\chi(i, j, y') \neq 0$.
    It follows that $\chi(i,j,y') = \chi(c,d,y') = -\chi(c', d', y)$.
\end{proposition}
Note that we do not assume that $\chi(i, j, y)$ is non-zero.
In particular for $i = c$ and $j = d$ we obtain $\chi(c',d',y) = -\chi(c, d, y')$.
\begin{proof}
    For $g \in \{c, c', i\}$ we inspect with Lemma~\ref{Lem:PointInSegment}
    \begin{align*}
        \chi(a,c,d)\chi(y',y,g) &= \ ?,\\
        \chi(y,c,d)\chi(a,y',g) &= 0,\\
        \chi(g,c,d)\chi(a,y,y') &= 0.
    \end{align*}
    As $\chi(y',c,d)\chi(a,y,g) \neq 0$ we conclude with Lemma~\ref{Lem:B2} that
    \[
        \chi(y',c,d)\chi(a,y,g) = \chi(a,c,d)\chi(y',y,g).
    \]
    With $\chi(a,y,c) = \chi(a,y,c') = \chi(a,y,i)$ we can therefore conclude that
    \[
        \chi(y',y,c) = \chi(y',y,c') = \chi(y',y,i).
    \]
    Analog
    \[
        \chi(y',y,d) = \chi(y',y,d') = \chi(y',y,j)
    \]
    and in particular by Lemma~\ref{Lem:PointInSegment}
    \[
        \chi(y',y,i) = \chi(y', y, c') = \chi(d', y, c') = -\chi(c', y, d') = - \chi(y', y, d') =   -\chi(y',y,j).
    \]
    Hence, $\chi(c', d', y) = -\chi(j, y', y)$.
    Now by assumptions $\chi(i,j,y) \neq -\chi(j, y', y)$.
    But as $\chi$ is acyclic one of $\chi(i, j, y)$, $\chi(j, y', y) = \chi(y', i, y) \neq 0$ needs to be equal to $\chi(i,j,y') \neq 0$.
    We conclude that
    \[
        \chi(i, j, y') = \chi(y', i, y) = \chi(j, y', y) = -\chi(c', d', y).
    \]
    Also
    \[
        \chi(y', i, y) = \chi(y', c, y) = \chi(y', c, d).
    \]
\end{proof}
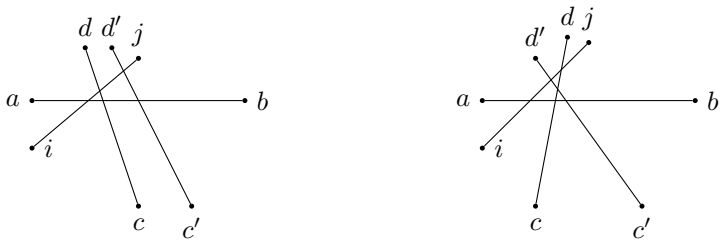
\begin{figure}
    \begin{center}
        \begin{tikzpicture}[scale=.7]
            \node (a) at (0, 0) [circle, fill=black, scale=0.2, label=left:$a$] {};
            \node (b) at ( 4, 0) [circle, fill=black, scale=0.2, label=right:$b$] {};
            \node (c) at (2, -2) [circle, fill=black, scale=0.2, label=below:$c$] {};
            \node (d) at (1,  1) [circle, fill=black, scale=0.2, label=above:$d$] {};
            \node (c1) at (3, -2) [circle, fill=black, scale=0.2, label=below:$c'$] {};
            \node (d1) at (1.5,  1) [circle, fill=black, scale=0.2, label=above:$d'$] {};
            \node (i) at (0, -0.9) [circle, fill=black, scale=0.2, label=right:$i$] {};
            \node (j) at (2, 0.8) [circle, fill=black, scale=0.2, label=above:$j$] {};
            \draw (a) -- (b);
            \draw (c) -- (d);
            \draw (c1) -- (d1);
            \draw (i) -- (j);
        \end{tikzpicture}
        \hspace{2cm}
        \begin{tikzpicture}[scale=.7]
            \node (a) at (0, 0) [circle, fill=black, scale=0.2, label=left:$a$] {};
            \node (b) at ( 4, 0) [circle, fill=black, scale=0.2, label=right:$b$] {};
            \node (c) at (1, -2) [circle, fill=black, scale=0.2, label=below:$c$] {};
            \node (d) at (1.6,  1.2) [circle, fill=black, scale=0.2, label=above:$d$] {};
            \node (c1) at (3, -2) [circle, fill=black, scale=0.2, label=below:$c'$] {};
            \node (d1) at (1,  0.8) [circle, fill=black, scale=0.2, label=above:$d'$] {};
            \node (i) at (0, -0.9) [circle, fill=black, scale=0.2, label=right:$i$] {};
            \node (j) at (2, 1.1) [circle, fill=black, scale=0.2, label=above:$j$] {};
            \draw (a) -- (b);
            \draw (c) -- (d);
            \draw (c1) -- (d1);
            \draw (i) -- (j);
        \end{tikzpicture}
    \end{center}
    \caption{Sketch of two cases of Proposition~\ref{Prop:Edges}.}
    \label{Fig:Edges}
\end{figure}
We will connect $v$ and $w$ unless we find a contradiction to Proposition~\ref{Prop:Edges}.

\smallskip

Given an acyclic chirotope $\restr{\chi}{X}$.
The small package
{\centering
\framebox{\url{https://github.com/kliem/TenColoredPoints}}\par
}
constructs the slightly enarged graph $H(\restr{\chi}{X})$ of $G(\restr{\chi}{X})$.
It then uses the algorithm from Section~\ref{Sec:kpkc} to verify that $H(\restr{\chi}{X})$ does not have a $k$-clique.
By Lemma~\ref{Lem:EnlargedGraph} this implies that $G(\restr{\chi}{X})$ also has no $k$-clique and this chirotope on $X$ satisfies the optimal colored Tverberg problem.

To iterate over all acyclic chirotopes on $10$ points of rank $3$ one can either use the list by~\cite{Aichholzer2001}\footnote{\url{http://www.ist.tugraz.at/staff/aichholzer/research/rp/triangulations/ordertypes/}.} or the package

{\centering
    \framebox{\url{https://github.com/kliem/pseudo_order_types}.}\par
}

Note that the list by~\cite{Aichholzer2001} does not contain the non-realizable acyclic chirotopes on $10$ points.
It would suffice to prove this instance of the optimal colored Tverberg problem,
but we also check the others to show it on the level of chirotopes.

\section{$k$-cliques in a $k$-partite graph}\label{Sec:kpkc}

\begin{table}
    \caption{
        Runtime in ms of checking for a $k$-clique for some graphs $H(\restr{\chi}{X})$.
    }
    \label{tab:SomeSamples}
    \begin{tabular}{r||rrrrr}
        Graph indexed by~\cite{Aichholzer2001} & \texttt{kpkc} & \texttt{FindClique} & \texttt{NetworkX} & \texttt{Cliquer} & \texttt{mcqd} \\
        \hline
        0                                      & 17,300        & nan                 & nan               & nan              & nan           \\
        1                                      & 17,300        & nan                 & nan               & nan              & nan           \\
        2                                      & 17,200        & nan                 & nan               & nan              & nan           \\
        20                                     & 5,070         & nan                 & nan               & nan              & nan           \\
        100                                    & 17,400        & nan                 & nan               & nan              & nan           \\
        1000                                   & 489           & 21,500              & nan               & nan              & nan           \\
        10000                                  & 178           & 3,620               & nan               & 626,000          & nan           \\
        1000000                                & 1,320         & nan                 & nan               & nan              & nan           \\
        2000000                                & 158           & 763                 & nan               & 9,440            & nan           \\
        5000000                                & 187           & 7,950               & nan               & nan              & nan           \\
        10000000                               & 22            & 6                   & nan               & 1,240            & nan           \\
    \end{tabular}
\end{table}

We provide a new algorithm to iterate over $k$-cliques in a $k$-partite graph.
It is implemented in \texttt{C++}:

{\centering
\framebox{\url{https://github.com/kliem/KPartiteKClique}.}\par
}

It is based on the depth-first algorithm of Grünert, Irnich, Zimmermann, Schneider and Wulfhorst~\cite{Grunert}:
\lstset{language    = Python,
        linewidth   = 0.95\textwidth,
        xleftmargin = 0.05\textwidth}
\begin{lstlisting}[belowskip=-0.8 \baselineskip]
def FINDCLIQUE(G):
    parts = G.parts()
    if len(parts) == 0:
        yield []
        return
    P = argmin(part.size() for part in parts)
    for v in P.vertices():
        V1 = v.neighbors()
        G1 = G.induced_subgraph(V1)
        for clique in FINDCLIQUE(G1):
            yield clique + [v]
\end{lstlisting}
This pivot selection is simple and fast for many purposes.
However, it did not terminate for some of the graphs encountered during this project.

Consider the graph

{\centering
\framebox{\url{https://github.com/kliem/PyKPartiteKClique/blob/main/sample_graphs/0.gz}.}\par
}
that can be recovered with \texttt{kpkc.test.load\_tester} from the Python wrapper of \texttt{KPartiteKClique}\footnote{
\url{https://github.com/kliem/PyKPartiteKClique}}.
This is $H(\restr{\chi}{X})$, where $\restr{\chi}{X}$ is the first order type in the enumeration by \cite{Aichholzer2001} listed on

{\centering
\framebox{\url{http://www.ist.tugraz.at/staff/aichholzer/research/rp/triangulations/ordertypes/}}\par
}
-- $10$ points in convex position.

The graph $H(\restr{\chi}{X})$ is a $71$-partite graph with 10,785~vertices and 6,630,275~edges.
The part coresponding to the colors has 10,045~vertices, the other parts have at most $20$~vertices.

Such a graph with density $0.11$ is unlikely to have a $71$-clique.
However, judging by the size of the graph, recursion depth $71$ is a challenge.
This graph has been tested with the following algorithms/implementations and none of them terminated in 24 hours:
\begin{itemize}
    \item \texttt{FindClique}~\cite{Grunert},
    \item \texttt{NetworkX}~\cite{Networkx1}~\cite{Networkx2},
    \item \texttt{Cliquer}~\cite{Cliquer},
    \item \texttt{mcqd}~\cite{Konc}.
\end{itemize}
Note that only \texttt{FindClique} exploits the given $71$-partition, while the other algorithms do note use it.
We provide our own implementation of \texttt{FindClique}, as the authors have not provided it in~\cite{Grunert}.

Proving Theorem~\ref{Thm:TenPoints} by Proposition~\ref{Prop:KClique} we need to check for $k$-cliques in 14,320,182 graphs.
Many of them are much simpler than this first one.
However, it seems desirable to find an algorithm that terminates for all those graphs in reasonable time.

To reduce the complexity of the graph, we start with vertices with few neighbors as done in \cite{Konc}:
\begin{lstlisting}[belowskip=-0.8 \baselineskip]
def kpkc(G, prec_depth=5):
    if len(parts) == 0:
        yield []
        return
    if prec_depth:
        G.sort_vertices(key=len_neighbors)
    for v in G.vertices():
        V1 = v.neighbors()
        G1 = G.induced_subgraph(V1)
        for clique in kpkc(G1, prec_depth - 1):
            yield clique + [v]
        G.remove(v)
\end{lstlisting}

The implementation is a bit more involved. In particular:
\begin{itemize}
    \item During sorting of the vertices we de facto remove vertices not connected to all parts.
        If this happend during the first sort, we sort again.
    \item After removing the second last vertex of a part, we sort again (and de facto select the last vertex).
    \item After removing the last vertex of a part, we immediatly return.
    \item The resources for the recursive calls are recycled to avoid memory allocations.
    \item The induced subgraph only keeps track of the selected vertices and their number of neighbors.
        The neighbors of a vertex are the intersection of the vertices of the subgraph with the neighbors of the original vertex.
        The number of neighbors is computed without storing the set of neighbors.
\end{itemize}

If a vertex has few neighbors, the induced subgraph of the neighbors is likely very easy to handle.
This reduces size of the graph vertex by vertex.
Table~\ref{tab:SomeSamples} compares the runtime for some graphs we need to analyse for Proposition~\ref{Prop:KClique}.
Apparantly, \texttt{kpkc} is the only choice suitable to solve our problem.

With \texttt{kpkc} we can verify in just 17 seconds that our first graph does not have a $71$-clique.
We prove Proposition~\ref{Prop:KClique} by analysing all 14,320,182 graphs.
Building and analysing all graphs was done in $780$ CPU-hours
using an Intel$^\text{\tiny{\textregistered}}$ Core\texttrademark{} i7-7700 CPU @ 3.60GHz x86\_64-processor.\footnote{
Actually, the 14,320,182 graphs were divided to multiple threads on multiple machines with the same specifications.}
This is an average of $196$ ms per instance.

Now, we will inspect benchmarks on random graphs and explain for what type of graphs the new algorithm is suitable.

\section{Benchmarks on random graphs}\label{Sec:Benchmarks}

In Tables~\ref{tab:SamplesType1a}, \ref{tab:SamplesType1b}, \ref{tab:SamplesType1c}, \ref{tab:SamplesType2a}, \ref{tab:SamplesType2b} and \ref{tab:SamplesType2c} we benchmark the implementations on random graphs.
Each row represents one graph randomly generated by certain parameters that is tested on all implementations.
A timeout after 1000 seconds is marked by ``nan''.
If an implementation is not included in the table, each entry would indicate a timeout.
We benchmark getting the first $k$-clique or checking for the existence of a such a $k$-clique.
We also benchmark obtaining all $k$-cliques, if the implementation has this available.
In Tables~\ref{tab:SamplesType2b} and~\ref{tab:SamplesType2c} we have not included timings for all cliques:
Either the graphs have few $k$-cliques and obtaing all $k$-cliques is just as fast as obtaining the first $k$-clique or the graph has many $k$-cliques and obtaining all of them results in a timout.

\subsection{Random graphs as constructed by~\cite{Grunert}}

\begin{table}
    \caption{
        Runtime in ms for sample graphs as constructed by \cite{Grunert}
    }
    \label{tab:SamplesType1a}
    \begin{tabular}{rrrrr||rrrrrrrr}
        $k$ & $\overset{\min}{|P_b|}$ & $\overset{\max}{|P_b|}$ & $a$  & $b$  & \multicolumn{2}{c}{\texttt{kpkc}} & \multicolumn{2}{c}{\texttt{FindClique}} & \multicolumn{2}{c}{\texttt{NetworkX}} & \texttt{Cliquer} & \texttt{mcqd} \\
            &                         &                         &      &      & first & all                       & first & all                             & first & all                           & first            & first         \\
        \hline
        5   & 50                      & 50                      & 0.14 & 0.14 & 1     & 1                         & 0     & 0                               & 11    & 11                            & 1                & 1             \\
        5   & 50                      & 50                      & 0.15 & 0.15 & 0     & 1                         & 0     & 0                               & 4     & 13                            & 1                & 1             \\
        5   & 50                      & 50                      & 0.2  & 0.2  & 0     & 1                         & 0     & 0                               & 2     & 24                            & 2                & 2             \\
        5   & 50                      & 50                      & 0.25 & 0.25 & 0     & 2                         & 0     & 1                               & 1     & 43                            & 2                & 2             \\
        5   & 50                      & 50                      & 0.0  & 0.3  & 0     & 1                         & 0     & 0                               & 1     & 13                            & 1                & 1             \\
        5   & 50                      & 50                      & 0.0  & 0.4  & 0     & 1                         & 0     & 0                               & 1     & 25                            & 2                & 2             \\
        5   & 50                      & 50                      & 0.0  & 0.45 & 0     & 3                         & 0     & 1                               & 1     & 43                            & 2                & 2             \\
        5   & 50                      & 50                      & 0.0  & 0.5  & 0     & 4                         & 0     & 2                               & 1     & 45                            & 2                & 2             \\
        \hline
        10  & 26                      & 37                      & 0.49 & 0.49 & 3     & 40                        & 0     & 1                               & 259   & 9,430                         & 46               & 45            \\
        10  & 26                      & 37                      & 0.5  & 0.5  & 0     & 35                        & 0     & 1                               & 35    & 8,260                         & 57               & 38            \\
        10  & 26                      & 37                      & 0.51 & 0.51 & 1     & 69                        & 0     & 1                               & 750   & 14,200                        & 110              & 66            \\
        10  & 26                      & 37                      & 0.4  & 0.6  & 2     & 47                        & 0     & 1                               & 217   & 13,900                        & 72               & 48            \\
        10  & 26                      & 37                      & 0.3  & 0.7  & 1     & 66                        & 0     & 4                               & 13    & 12,700                        & 152              & 51            \\
        10  & 50                      & 50                      & 0.42 & 0.42 & 13    & 93                        & 0     & 1                               & 1,090 & 21,400                        & 107              & 113           \\
        10  & 50                      & 50                      & 0.43 & 0.43 & 4     & 121                       & 0     & 2                               & 1,700 & 27,800                        & 180              & 131           \\
        10  & 50                      & 50                      & 0.44 & 0.44 & 22    & 151                       & 0     & 2                               & 316   & 34,600                        & 214              & 170           \\
        10  & 50                      & 50                      & 0.46 & 0.46 & 5     & 304                       & 0     & 4                               & 27    & 55,200                        & 324              & 275           \\
        10  & 50                      & 50                      & 0.48 & 0.48 & 1     & 573                       & 0     & 8                               & 49    & 83,900                        & 404              & 395           \\
        10  & 50                      & 50                      & 0.5  & 0.5  & 1     & 997                       & 0     & 20                              & 16    & 172,000                       & 1,160            & 694           \\
    \end{tabular}
\end{table}
\begin{table}
    \caption{
        Runtime in ms for sample graphs as constructed by \cite{Grunert}
    }
    \label{tab:SamplesType1b}
    \begin{tabular}{rrrrr||rrrrrrrr}
        $k$ & $\overset{\min}{|P_b|}$ & $\overset{\max}{|P_b|}$ & $a$   & $b$   & \multicolumn{2}{c}{\texttt{kpkc}} & \multicolumn{2}{c}{\texttt{FindClique}} & \texttt{Cliquer} & \texttt{mcqd} \\
            &                         &                         &       &       & first   & all                     & first & all                             & first            & first         \\
        \hline
        50  & 5                       & 15                      & 0.91  & 0.91  & nan     & nan                     & 64    & 64                              & nan              & nan           \\
        50  & 5                       & 15                      & 0.918 & 0.918 & nan     & nan                     & 28    & 10,600                          & nan              & nan           \\
        50  & 5                       & 15                      & 0.92  & 0.92  & nan     & nan                     & 8     & 2,910                           & nan              & nan           \\
        \hline
        20  & 23                      & 39                      & 0.7   & 0.7   & 187,000 & 320,000                 & 106   & 193                             & nan              & nan           \\
        20  & 23                      & 39                      & 0.71  & 0.71  & 40,700  & 535,000                 & 52    & 345                             & nan              & nan           \\
        20  & 23                      & 39                      & 0.72  & 0.72  & 5,120   & 955,000                 & 2     & 642                             & nan              & nan           \\
        20  & 23                      & 39                      & 0.7   & 0.73  & 2,570   & nan                     & 4     & 867                             & nan              & nan           \\
        20  & 23                      & 39                      & 0.65  & 0.78  & 1,040   & 454,000                 & 0     & 436                             & nan              & nan           \\
        \hline
        30  & 11                      & 30                      & 0.6   & 0.6   & 429     & 429                     & 0     & 0                               & 796,000          & 121,000       \\
        30  & 11                      & 30                      & 0.7   & 0.7   & 16,300  & 16,300                  & 1     & 1                               & nan              & nan           \\
        30  & 11                      & 30                      & 0.8   & 0.8   & nan     & nan                     & 1,160 & 1,160                           & nan              & nan           \\
        30  & 11                      & 30                      & 0.81  & 0.81  & nan     & nan                     & 1,330 & 1,330                           & nan              & nan           \\
        30  & 11                      & 30                      & 0.82  & 0.82  & nan     & nan                     & 436   & 3,450                           & nan              & nan           \\
        30  & 11                      & 30                      & 0.84  & 0.84  & nan     & nan                     & 6     & 51,300                          & nan              & nan           \\
        30  & 11                      & 30                      & 0.88  & 0.88  & 1,440   & nan                     & 0     & nan                             & nan              & nan           \\
        \hline
        100 & 10                      & 10                      & 0.7   & 0.7   & 52      & 52                      & 0     & 0                               & nan              & nan           \\
        100 & 10                      & 10                      & 0.8   & 0.8   & 4,760   & 4,760                   & 0     & 0                               & nan              & nan           \\
        100 & 10                      & 10                      & 0.85  & 0.85  & 221,000 & 221,000                 & 2     & 2                               & nan              & nan           \\
        100 & 10                      & 10                      & 0.9   & 0.9   & nan     & nan                     & 149   & 149                             & nan              & nan           \\
        100 & 10                      & 10                      & 0.92  & 0.92  & nan     & nan                     & 4,190 & 4,190                           & nan              & nan           \\
        100 & 10                      & 10                      & 0.94  & 0.94  & nan     & nan                     & nan   & nan                             & nan              & nan           \\
        100 & 10                      & 10                      & 0.95  & 0.95  & nan     & nan                     & nan   & nan                             & nan              & nan           \\
        100 & 10                      & 10                      & 0.97  & 0.97  & nan     & nan                     & 3     & nan                             & nan              & nan           \\
    \end{tabular}
\end{table}
\begin{table}
    \caption{
        Runtime in ms for sample graphs as constructed by \cite{Mirghorbani}
    }
    \label{tab:SamplesType1c}
    \begin{tabular}{rrrrr||rrrrrrrr}
        $k$ & $\overset{\min}{|P_b|}$ & $\overset{\max}{|P_b|}$ & $a$  & $b$  & \multicolumn{2}{c}{\texttt{kpkc}} & \multicolumn{2}{c}{\texttt{FindClique}} & \multicolumn{2}{c}{\texttt{NetworkX}} & \texttt{Cliquer} & \texttt{mcqd} \\
            &                         &                         &      &      & first & all                       & first & all                             & first & all                           & first            & first         \\
        \hline
        3   & 100                     & 100                     & 0.1  & 0.1  & 0     & 3                         & 0     & 2                               & 1     & 7                             & 1                & 1             \\
        4   & 100                     & 100                     & 0.15 & 0.15 & 0     & 5                         & 0     & 2                               & 2     & 40                            & 4                & 3             \\
        5   & 100                     & 100                     & 0.2  & 0.2  & 0     & 12                        & 0     & 2                               & 5     & 190                           & 9                & 9             \\
        6   & 100                     & 100                     & 0.25 & 0.25 & 0     & 36                        & 0     & 3                               & 8     & 996                           & 23               & 23            \\
        \hline
        7   & 50                      & 50                      & 0.35 & 0.35 & 0     & 16                        & 0     & 1                               & 8     & 738                           & 13               & 12            \\
        8   & 50                      & 50                      & 0.4  & 0.4  & 0     & 44                        & 0     & 2                               & 12    & 3,890                         & 38               & 34            \\
        9   & 50                      & 50                      & 0.45 & 0.45 & 0     & 174                       & 0     & 4                               & 9     & 20,500                        & 171              & 125           \\
        10  & 50                      & 50                      & 0.5  & 0.5  & 1     & 925                       & 0     & 17                              & 143   & 131,000                       & 1,480            & 643           \\
    \end{tabular}
\end{table}

We first rerun the benchmarks by \cite[Table~2]{Grunert}: A random $k$-partite graph is generated by parameters
\[
    (k, \min |P_b|, \max |P_b|, a, b).
\]
Each part has a random number of vertices in $\{\min |P_b|, \min |P_b| + 1, \dots, \max |P_b|\}$ by uniform distribution.
To each vertex $v$ we associate a random number $p_v$, which is uniformly selected from the interval $[a, b]$.
Finally, two vertices $v, w$ are connected by an edge with probability $\frac{p_v + p_w}{2}$.

Indeed,~\texttt{FindClique} is by far the best choice for those graphs as can be seen in Tables~\ref{tab:SamplesType1a} and~\ref{tab:SamplesType1b}.
While~\cite{Grunert} only obtains the first 1000 cliques, we try to find all cliques.
They used a 100 MHz machine with 32 MB of RAM.
Thus already the advance of technology has improved the benchmarks by a factor of at least 36.
However, for increasing $k$ our implementation appears to be an improvement.
For $k = 100$ we have improved the old timings by a factor of 1000, which is 28-times faster than that factor of 36.

In an intermediate paper Mirghorbani and Krokhmal~\cite{Mirghorbani} proposed to improve the data structure of \texttt{FindClique}.
They reported that using arrays and bitsets could each gain a factor of about $3$.
In~\ref{tab:SamplesType1c} we have rerun the tests of~\cite[Table~1]{Mirghorbani}, which suggests that our implementation is yet faster by a factor of a bit more than $3$
when considering that they only used a 3 GHz machine.
If~\cite{Mirghorbani} improved \texttt{FindClique} of up to 9 and our implementation is yet an improvement of up to 3, this agrees with the above observation for $k=100$.

However, the main advantage is that we now have a published implementation of~\texttt{FindClique},
which neither~\cite{Mirghorbani} nor~\cite{Grunert} have provided.

\subsection{Rare attraction random graphs}

\begin{table}
    \caption{
        Runtime in ms for rare attraction random graphs
    }
    \begin{tabular}{rrr||rrrrrrrr}
        $k$ & $\overset{\max}{|P_b|}$ & $a$  & \multicolumn{2}{c}{\texttt{kpkc}} & \multicolumn{2}{c}{\texttt{FindClique}} & \multicolumn{2}{c}{\texttt{NetworkX}} & \texttt{Cliquer} & \texttt{mcqd} \\
            &                         &      & first & all                       & first & all                             & first & all                           & first            & first         \\
        \hline
        5   & 10                      & 0.1  & 0     & 0                         & 0     & 0                               & 0     & 0                             & 0                & 0             \\
        5   & 10                      & 0.2  & 0     & 0                         & 0     & 0                               & 0     & 1                             & 0                & 0             \\
        5   & 20                      & 0.05 & 0     & 1                         & 0     & 1                               & 0     & 4                             & 0                & 0             \\
        5   & 20                      & 0.1  & 0     & 3                         & 0     & 3                               & 0     & 6                             & 0                & 0             \\
        5   & 50                      & 0.01 & 0     & 92                        & 0     & 77                              & 1     & 155                           & 1                & 1             \\
        5   & 50                      & 0.02 & 0     & 111                       & 0     & 94                              & 1     & 166                           & 1                & 1             \\
        \hline
        10  & 10                      & 0.4  & 0     & 4                         & 0     & 2                               & 0     & 22                            & 0                & 0             \\
        10  & 10                      & 0.6  & 0     & 22                        & 0     & 12                              & 0     & 51                            & 0                & 0             \\
        10  & 20                      & 0.3  & 0     & 50                        & 0     & 27                              & 2     & 790                           & 1                & 1             \\
        10  & 20                      & 0.5  & 0     & 1,160                     & 0     & 780                             & 1     & 3,280                         & 1                & 3             \\
        10  & 50                      & 0.05 & 0     & 49                        & 0     & 66                              & 19    & 36,800                        & 4                & 5             \\
        10  & 50                      & 0.1  & 0     & 237                       & 0     & 189                             & 5     & 47,600                        & 4                & 5             \\
        10  & 100                     & 0.01 & 0     & 5,870                     & 0     & 7,790                           & 19    & nan                           & 17               & 73            \\
        10  & 100                     & 0.02 & 0     & 8,490                     & 0     & 9,930                           & 14    & nan                           & 18               & 76            \\
    \end{tabular}
    \label{tab:SamplesType2a}
\end{table}
\begin{table}
    \caption{
        Runtime in ms for rare attraction random graphs
        to check for a $k$-clique
    }
    \label{tab:SamplesType2b}
    \begin{tabular}{rrr||rrrrr}
        $k$ & $\overset{\max}{|P_b|}$ & $a$  & \texttt{kpkc} & \texttt{FindClique} & \texttt{Cliquer} \\
        \hline
        50  & 20                      & 0.5  & 208           & 20                  & nan              \\
        50  & 20                      & 0.6  & 2,520         & 422                 & nan              \\
        50  & 20                      & 0.7  & 158,000       & 23,800              & nan              \\
        50  & 20                      & 0.71 & 304,000       & 248,000             & nan              \\
        50  & 20                      & 0.72 & nan           & 712,000             & nan              \\
        50  & 20                      & 0.73 & nan           & 761,000             & nan              \\
        50  & 20                      & 0.75 & nan           & nan                 & nan              \\
        50  & 20                      & 0.76 & nan           & nan                 & nan              \\
        50  & 20                      & 0.77 & nan           & 15,600              & nan              \\
        50  & 20                      & 0.78 & nan           & 1,500               & nan              \\
        50  & 20                      & 0.79 & nan           & 189                 & nan              \\
        50  & 20                      & 0.8  & nan           & 6                   & nan              \\
        \hline
        50  & 50                      & 0.1  & 36            & 4,460               & nan              \\
        50  & 50                      & 0.2  & 128           & 33,300              & nan              \\
        50  & 50                      & 0.3  & 2,090         & 174,000             & nan              \\
        50  & 50                      & 0.4  & 24,600        & nan                 & nan              \\
        50  & 50                      & 0.5  & 803,000       & nan                 & nan              \\
        50  & 50                      & 0.6  & nan           & nan                 & nan              \\
        50  & 50                      & 0.71 & nan           & nan                 & nan              \\
        50  & 50                      & 0.72 & nan           & 29,900              & 183,000          \\
        50  & 50                      & 0.73 & nan           & 11,700              & nan              \\
        50  & 50                      & 0.74 & nan           & 1,550               & nan              \\
        \hline
        50  & 100                     & 0.1  & 243           & nan                 & nan              \\
        50  & 100                     & 0.2  & 11,100        & nan                 & nan              \\
        50  & 100                     & 0.3  & 235,000       & nan                 & nan              \\
        50  & 100                     & 0.4  & nan           & nan                 & nan              \\
        50  & 100                     & 0.64 & nan           & nan                 & nan              \\
        50  & 100                     & 0.65 & nan           & nan                 & 119,000          \\
        50  & 100                     & 0.66 & nan           & nan                 & nan              \\
        50  & 100                     & 0.67 & nan           & nan                 & nan              \\
        50  & 100                     & 0.68 & nan           & nan                 & nan              \\
        50  & 100                     & 0.69 & nan           & 22,800              & nan              \\
        50  & 100                     & 0.7  & nan           & 2,270               & 42,900           \\
    \end{tabular}
\end{table}
\begin{table}
    \caption{
        Runtime in ms for rare attraction random graphs
        to check for a $k$-clique
    }
    \label{tab:SamplesType2c}
    \begin{tabular}{rrr||rrrrr}
        $k$ & $\overset{\max}{|P_b|}$ & $a$  & \texttt{kpkc} & \texttt{FindClique} \\
        \hline
        100 & 20                      & 0.4  & 48            & 2                   \\
        100 & 20                      & 0.5  & 435           & 13                  \\
        100 & 20                      & 0.6  & 4,650         & 98                  \\
        100 & 20                      & 0.7  & 475,000       & 9,960               \\
        100 & 20                      & 0.8  & nan           & nan                 \\
        100 & 20                      & 0.89 & nan           & nan                 \\
        100 & 20                      & 0.9  & nan           & 2,430               \\
        \hline
        100 & 50                      & 0.1  & 178           & 7,150               \\
        100 & 50                      & 0.2  & 292           & 81,900              \\
        100 & 50                      & 0.3  & 8,120         & 496,000             \\
        100 & 50                      & 0.4  & 46,100        & nan                 \\
        100 & 50                      & 0.5  & nan           & nan                 \\
        100 & 50                      & 0.87 & nan           & nan                 \\
        100 & 50                      & 0.88 & nan           & 4,240               \\
        100 & 50                      & 0.89 & nan           & 2                   \\
        100 & 50                      & 0.9  & nan           & 0                   \\
        \hline
        100 & 100                     & 0.1  & 892           & nan                 \\
        100 & 100                     & 0.2  & 41,600        & nan                 \\
        100 & 100                     & 0.3  & 260,000       & nan                 \\
        100 & 100                     & 0.4  & nan           & nan                 \\
        100 & 100                     & 0.85 & nan           & nan                 \\
        100 & 100                     & 0.86 & nan           & 7,750               \\
        100 & 100                     & 0.87 & nan           & 4                   \\
        100 & 100                     & 0.88 & nan           & 0                   \\
        100 & 100                     & 0.89 & nan           & 0                   \\
        100 & 100                     & 0.9  & nan           & 0                   \\
    \end{tabular}
\end{table}

It seems that \texttt{FindClique} is the algorithm of choice for all random graphs as constructed by~\cite{Grunert}.
The reason seems to be that edges are somewhat equally distributed.
Vertices in small parts have about the same expected number of neighbors.
This makes the pivot selection of \texttt{kpkc} useless.
However, $k$-partite graphs that correspond to real life problems might behave differently:

Suppose there is only one cement mill in the area, two concrete pumps, twenty conrete mixer trucks, and twenty concrete crews.
Nobody can question the quality of the cement mill, because there is no alternative.
As there is only two concrete pumps, the truck drivers will usually be willing to work with both of them.
Likewise the concrete crews will usually put up with both pump operators.
However, it is very much possible that the conrete crews might refuse to work with some truck drivers (always late) or the truck drivers might refuse to work with some crews (always order more trucks than they need).

This problem corresponds to a $4$-partite graph.
\texttt{FindClique} first selects the cement mill, which is trivial.
As a next step it divides the problem in two:
$4$-cliques containing one pump and $4$-cliques containing the other.
However, this is only a good choice if the pumps have very different sets of neighbors.
The real problem is assigning truck drivers to concrete crews, which \texttt{FindClique} solves twice, once for each pump.
On the other hand, \texttt{kpkc} immediately solves the actual problem..

The graphs $H(\restr{\chi}{X})$ seem to behave somewhat like this example.
The less choices remain for a valid intersection point $y$, the less information such a selection will gain.
Thus, \texttt{FindClique} is likely to divide the problem into two or more very similar problems.
\texttt{kpkc} instead selects choices that are unlikely to correspond to counter examples and quickly rules them out.

We construct \defn{rare attraction random graphs} parametrized by $(k, \max |P_b|, a)$ as follows:
The part $i$ has size
\[
    1 + \lfloor i\frac{(\max |P_b| - 1)}{k} \rfloor
\]
for $i = 1,\dots,k$.
Let $f$ be the affine function determined by $f(1) = 1$ and $f(\max |P_b|) = a$.
We generate an edge between vertices of different parts of sizes $s$ and $t$ with probability $f(\min(s, t))$.
The fewer vertices a part has, the rarer those vertices are.
If vertices are rare, other vertices are more likely to be attracted.

Table~\ref{tab:SamplesType2a} reveals that for smaller graphs \texttt{FindClique} is still the algorithm of choice.
In Tables~\ref{tab:SamplesType2b} and~\ref{tab:SamplesType2c} we see that for $k \geq 50$, $\max |P_b| \geq 50$ and low density \texttt{kpkc} is much better in finding $k$-cliques (or rather verifying their abscence).
With high density and probably many $k$-cliques, \texttt{FindClique} is still much faster (at least in finding some $k$-cliques).

\subsection{Conclusion}

For finding $k$-cliques in $k$-partite graphs it is recommended to use a specialized algorithm.
If the graph is large and expected to have few $k$-cliques and parts with fewer vertices have more neighbors, then \texttt{kpkc} is probably the better algorithm.
In most other cases, \texttt{FindClique} seems better suited.
In either case, an algorithm without implementation requires lots of work to be of any use.
With

{\centering
\framebox{\url{https://github.com/kliem/KPartiteKClique}}\par
}
both algorithms are implemented using static polymorphism.

\section*{Acknowledgements}
I would like to thank Jean-Philippe Labbé for much help during the start of this project and in particular for pointing out how cliques can be useful.
Further, I would like to thank Maikel Nadolski for many valuable discussions and much help regarding the implementation in \texttt{C++}.
Additionally, I would like to thank Florian Frick for many valuable comments and for providing references.

\end{document}